\magnification=1200

\def\C{{\bf C}}
\def\K{{\bf K}}
\def\N{{\bf N}}
\def\R{{\bf R}}
\def\span{{\rm span}}
\def\hal{{\vrule height 10pt width 4pt depth 0pt}}

\centerline{\bf Set theory and cyclic vectors}
\medskip

\centerline{Nik Weaver\footnote{*}{Supported by NSF grant DMS-0070634\hfill\break
Math Subject Classification numbers: 03E35, 03E50, 47A15, 47A16}}
\bigskip
\bigskip

{\narrower{
\noindent \it Let $H$ be a separable, infinite dimensional Hilbert space
and let $S$ be a countable subset of $H$. Then most positive operators on
$H$ have the property that every nonzero vector in the span of $S$ is cyclic,
in the sense that the set of operators in the positive part of the unit ball
of $B(H)$ with this property is comeager for the strong operator topology.
\medskip

Suppose $\kappa$ is a regular cardinal such that $\kappa \geq \omega_1$
and $2^{<\kappa} = \kappa$. Then it is relatively consistent with ZFC that
$2^\omega = \kappa$ and for any subset $S \subset H$ of cardinality less
than $\kappa$ the set of positive operators in the unit ball of $B(H)$ for
which every nonzero vector in the span of $S$ is cyclic is comeager for the
strong operator topology.
\bigskip}}
\bigskip

\noindent {\bf 1. Introduction}
\bigskip

Let $H$ be a separable, infinite dimensional Hilbert space and let $B(H)$
be the set of bounded linear operators $A: H \to H$. A closed subspace $E$
of $H$ is {\it invariant} for such an operator $A$ if $A(E) \subset E$. The
invariant subspace problem (ISP) for Hilbert spaces asks whether there
exists an operator $A \in B(H)$ whose only closed invariant subspaces are
$\{0\}$ and $H$.
\medskip

It was shown by Enflo [2] that there exists a bounded operator on a Banach
space that has no proper closed invariant subspaces. Read [7] showed that
the Banach space could be taken to be $l^1 = l^1(\N)$. A simplified version
of Read's example is given in [1]. The ISP remains open for Hilbert spaces;
it is also unknown whether there exists a separable, infinite dimensional
Banach space on which every bounded operator has a proper closed invariant
subspace.
\medskip

This note was motivated by the observation that the ISP for Hilbert spaces
can be rephrased as a question about the existence of a generic filter on
a certain poset. (This material is not needed for Section 2.) The construction
is this. Let $P$ be the poset consisting of all partially defined operators
$A$ on the Hilbert space $l^2 = l^2(\N)$ with the properties
\medskip

(a) ${\rm dom}(A)$ is a finite dimensional subspace of $l^2$;

(b) if $E$ is a subspace of ${\rm dom}(A)$ and $A(E) \subset E$ then
$E = \{0\}$; and

(c) $\|A\| < 1$.
\medskip

\noindent Order $P$ by reverse inclusion. For any unit vectors $v, w
\in l^2$ define
$$D_{v,w} = \{A \in P: \hbox{ there exists $n$ such that $A^n(v)$ is
defined and }\langle A^n v, w\rangle \neq 0\}.$$
It is not too hard to see that every $D_{v,w}$ is dense in $P$, and a
filter of $P$ which intersects every $D_{v,w}$ defines a bounded
operator with no proper invariant subspaces. Conversely, if there is
such an operator it can be scaled to have norm $< 1$, and then its
finite dimensional restrictions define a filter of $P$ which intersects
every $D_{v,w}$. Thus, the ISP for Hilbert space can be cast in
set-theoretic terms: it is equivalent to the existence of a $D$-generic
filter on $P$, where $D = \{D_{v,w}: v, w \in l^2$ and $\|v\| = \|w\| = 1\}$.
\medskip

The poset $P$ is not ccc, but this is not essential; for example, it can
be replaced by the countable poset of all finite matrices with rational
entries, ordered by a reasonable notion of approximate inclusion. Thus,
one can apply Martin's axiom (see, e.g., [6]) to obtain the consistency
of an operator which meets ``many'' of the sets $D_{v,w}$. This raises the
possibility that the ISP for Hilbert space may be independent of ZFC.
However, the assertion that $A \in B(l^2)$ has no invariant subspaces is
$\Sigma_1$ (it can be reformulated as ``for all unit vectors $v,w \in l^2$
there exists $n$ such that $\langle A^n v,w\rangle \neq 0$'') and hence
absolute, so if the ISP is independent, this cannot be shown using a
straightforward forcing argument.
\bigskip
\bigskip

\noindent {\bf 2. A relative consistency result}
\bigskip

As we indicated above, although Martin's axiom alone will not suffice to
prove the consistency of an operator with no proper closed invariant subspaces
(unless this can already be proven in ZFC), it does allow one to prove partial
results in this direction. In this section we present perhaps the strongest
natural result along these lines. It was in fact originally proven directly
from Martin's axiom, but here we give a better proof based on a suggestion
of Kenneth Kunen.
\medskip

It is easy to see that the operator $A \in B(H)$ has no proper closed
invariant subspaces if and only if every nonzero vector is cyclic, i.e.,
for every nonzero $v \in H$ the span of the sequence $(A^nv)$ is dense
in $H$. Thus, the more cyclic vectors $A$ has, the ``closer'' it gets to
being a counterexample to the ISP.
\medskip

If $X$ is a Banach space then we let $[X]_1$ denote its closed unit ball. The
{\it strong operator topology} on $B(H)$ is generated by the basic open sets
$${\cal O}_{B, E, \epsilon} = \{A \in B(H): \|A|_E - B\| < \epsilon\}$$
taken over $\epsilon > 0$, $E$ a finite dimensional subspace of $H$, and
$B: E \to H$ a linear map. If $H$ is separable then this topology is
second countable and its restriction to $[B(H)]_1$ is metrizable by
$d(A,A') = \sum 2^{-n} \|(A - A')(v_n)\|$, where $(v_n)$ is a countable
dense subset of $[H]_1$. Moreover, this metric is complete, so $[B(H)]_1$
with the relative strong operator topology is a Polish space.
\medskip

Let $[B(H)]^+_1$ denote the set of positive operators $A \in [B(H)]_1$,
i.e., those self-adjoint operators which satisfy $0 \leq \langle Av, v\rangle
\leq 1$ for all $v \in H$. (In the case of complex scalars, self-adjointness
follows from the second condition.) If $H$ is separable then $[B(H)]^+_1$ is
a Polish space via the same metric which shows that $[B(H)]_1$ is Polish.
\bigskip

\noindent {\bf Lemma.} {\it Let $H$ be a separable, infinite dimensional
Hilbert space and let $E$ be a finite dimensional subspace of $H$.
Then the set of operators in $[B(H)]^+_1$ for which every nonzero vector in
$E$ is cyclic is comeager for the relative strong operator topology.}
\medskip

\noindent {\it Proof.} Let $(x_n)$ be an orthonormal basis of $H$. For
$m \in \N$ and $\delta > 0$ let $U_{m,\delta}$ be the set of operators
$A \in [B(H)]^+_1$ such that
$$d(x_m, \span\{A^k v: k \in \N\}) < \delta$$
for every nonzero $v \in E$. We will show that $U_{m,\delta}$ is open and
dense in $[B(H)]^+_1$ for every $m$ and $\delta$. Intersecting the sets
$U_{m,\delta}$ over all $m$ and all $\delta$
of the form $\delta = 1/n$ yields the set of $A \in [B(H)]^+_1$ for
which $\span\{A^k v: k \in \N\} = H$ for every nonzero $v \in E$;
so we will have shown that this set is a countable intersection of
open, dense sets, as desired.
\medskip

Fix $m$ and $\delta$ for the remainder of the proof. We first show
that $U_{m,\delta}$ is open. Let $A \in U_{m,\delta}$ and for each
unit vector $v \in E$ let $f(v)$ be the smallest integer such that
$$d(x_m, \span\{A^k v: 0 \leq k \leq f(v)\}) < \delta.$$
Then the function $f$ is upper semicontinuous on the unit sphere of $E$
(which is compact), so $f$ is bounded. Let $N$ be an upper bound for $f$
and let $F = \span\{A^k v: v \in E, 0 \leq k \leq N\}$ and
$$\delta' = \sup \{d(x_m, \span\{A^k v: 0 \leq k \leq N\}):
v \in E, \|v\| = 1\}.$$
Since $E$ is finite dimensional, so is $F$. Also, by compactness of the
unit sphere $\delta' < \delta$. Now for $\epsilon > 0$ let $U_\epsilon$
be the set of operators $B \in [B(H)]^+_1$ which satisfy $\|(B - A)w\|
< \epsilon$ for all $w \in [F]_1$. This set is strong operator open for
each positive $\epsilon$. For each unit vector $v \in E$ let $g(v)$ be the
supremum of the set of $\epsilon > 0$ such that $B \in U_\epsilon$ implies
$$d(x_m, \span\{B^k v: 0 \leq k \leq N\})
\leq {1\over 2}(\delta' + \delta) < \delta.$$
Then every vector in the unit sphere of $F$ has a neighborhood in which
$g$ is bounded away from $0$, so $g$ must be bounded below by some positive
$\epsilon$, and we have $A \in U_\epsilon \subset U_{m,\delta}$ for this
$\epsilon$. Thus $U_{m,\delta}$ is strong operator open.
\medskip

Now we must show that $U_{m,\delta}$ is strong operator dense in $[B(H)]^+_1$.
Fix $A \in [B(H)]^+_1$, a finite dimensional subspace $F$ of $H$ which
contains $E$ and $x_m$, and $\epsilon > 0$. We will find an operator
$B \in U_{m,\delta}$ such that $\|(B - A)w\| < \epsilon$ for all
$w \in [F]_1$. Let $F' = \span(F + A(F))$, let $n = {\rm dim}(F')$,
let $P_{F'}$ be the orthogonal projection of $H$ onto $F'$, and let
$A' = P_{F'} A P_{F'}$. Note that $A'v = Av$ for all $v \in F$. Next,
choose an integer $r > 4n/\delta^2$ and let $X$ be a subspace of $H$ of
dimension $nr$ which contains $F'$. We can identify $X$ with the tensor
product space $\K^n \otimes \K^r$ (where $\K$ is the scalar field, $\K = \R$
or $\K = \C$) in such a way that $F'$ is identified with
$\K^n \otimes \{(r^{-1/2}, \ldots, r^{-1/2})\}$.
\medskip

Since $A$ is positive and $\|A\| \leq 1$, the same is true of $A'$.
Thus $F'$ is spanned by eigenvectors of $A'$, each of which
belongs to an eigenvalue between 0 and 1, inclusive. Let
$v_1', \ldots, v_n'$ be an orthonormal set of eigenvectors
belonging to the eigenvalues $\lambda_1, \ldots, \lambda_n$. (The
$\lambda_i$ need not be distinct.) According to the above identification
we have $v_i' = v_i \otimes (r^{-1/2}, \ldots, r^{-1/2})$ for some
orthonormal basis $\{v_i\}$ of $\K^n$.
\medskip

Let $\{w_j: 1 \leq j \leq r\}$ be the standard orthonormal basis of $\K^r$.
Then the vectors $v_i \otimes w_j$ constitute an orthonormal basis of
$X \cong \K^n \otimes \K^r$. Define $B' \in B(X)$ by setting
$B'(v_i\otimes w_j) = \lambda_i v_i\otimes w_j$. Then $B'$ is positive,
$\|B'\| \leq 1$, and $B'|_F = A'|_F = A|_F$.
\medskip

To complete the proof, we will show that there exists $B \in [B(X)]^+_1$
such that $\|B - B'\| < \epsilon$ and $B P_X \in U_{m, \delta}$. We will
define $B$ by choosing an orthonormal basis $\{e_{ij}\}$ of $X$ and
corresponding values $0 \leq \sigma_{ij} \leq 1$ and setting $Be_{ij}
= \sigma_{ij}e_{ij}$. If each $e_{ij}$ is sufficiently close to
$v_i\otimes w_j$ and each $\sigma_{ij}$ is sufficiently close to
$\lambda_i$ then we will have $\|B - B'\| < \epsilon$. Thus, we must
show that there exist $\{e_{ij}\}$ and $\{\sigma_{ij}\}$ arbitrarily
close to $\{v_i\otimes w_j\}$ and $\{\lambda_i\}$ which achieve
$B P_X \in U_{m, \delta}$.
\medskip

First, we claim that there exist orthonormal bases $\{e_{ij}\}$
arbitrarily close to the basis $\{v_i\otimes w_j\}$ with the
property that every $n$-element subset of the set $\{P_{F'}(e_{ij})\}$
is linearly independent. That is, any $n$ vectors in the basis
orthogonally project to independent vectors in $F'$. This is true
because, for any $n$ indices $i_1j_1, \ldots, i_nj_n$ the family
of bases $\{e_{ij}\}$ for which the vectors $P_{F'}(e_{i_1j_1}),
\ldots, P_{F'}(e_{i_nj_n})$ are dependent is a variety of codimension 1
in the manifold of all orthonormal bases. Thus, the family of bases
for which some $n$ elements project onto a dependent set is a
union of $rn \choose n$ meager sets, and hence meager. So, we can
perturb the basis $\{v_i \otimes w_j\}$ by an arbitrarily small
amount and achieve this condition.
\medskip

Now, having chosen $\{e_{ij}\}$ so as to satisfy the previous
claim, we conclude by showing that any choice of {\it distinct} values
$\sigma_{ij}$ such that each difference $|\sigma_{ij} - \lambda_i|$
is sufficiently small will ensure $B P_X \in U_{m,\delta}$.
To see this, observe first that for any nonzero $v \in F'$ at most
$n - 1$ of the inner products $\langle v, e_{ij}\rangle$ are zero.
Otherwise, $n$ of the vectors $e_{ij}$ would be orthogonal to
$v$, and hence $n$ of the vectors $P_{F'}(e_{ij})$ would be orthogonal
to $v$, which would imply linear dependence since ${\rm dim}(F') = n$.
This contradicts the previous claim. Now for each nonzero $v \in F'$ let
$$F_v = \span\{e_{ij}: \langle v, e_{ij}\rangle \neq 0\}.$$
Distinctness of the $\sigma_{ij}$ implies that the vectors
$B^k v$ are linearly independent for $0 \leq k < {\rm dim}(F_v)$;
since $F_v$ clearly contains $\span\{B^k v: k \in \N\}$ (it contains
$v$ and is invariant for $B$) this shows that the two are equal. Thus,
we must show that $d(x_m, F_v) < \delta$. But $x_m \in F'$, so
$|\langle x_m, v_i \otimes w_j\rangle| \leq r^{-1/2}$
for every $i$ and $j$. We may therefore assume that
$|\langle x_m, e_{ij}\rangle| \leq 2r^{-1/2} < \delta/\sqrt{n}$ for every
$i$ and $j$. Thus, if $x_m' = x_m - P_{F_v}(x_m)$ then
$$\|x_m'\|^2 = \sum_{e_{ij} \not\in F_v} |\langle x_m, e_{ij}\rangle|^2
< n\cdot (\delta^2/n) = \delta^2,$$
i.e., $\|x_m'\| < \delta$, since $F_v$ contains all but at most $n - 1$ of
the vectors $e_{ij}$. This shows that $d(x_m, F_v) < \delta$, as
desired.\hfill\hal
\bigskip

Since a countable intersection of comeager sets is comeager, the lemma
immediately implies the following result.
\bigskip

\noindent {\bf Theorem 1.} {\it Let $H$ be a separable, infinite dimensional
Hilbert space and let $S$ be a countable subset of $H$. Then the set of
operators in $[B(H)]^+_1$ for which every nonzero vector in the span of $S$
is cyclic is comeager in $[B(H)]^+_1$.}
\bigskip

\noindent {\bf Theorem 2.} {\it Let $H$ be a separable, infinite dimensional
Hilbert space and let $\kappa$ be a regular cardinal such that $\kappa \geq
\omega_1$ and $2^{<\kappa} = \kappa$. Then it is relatively consistent with
ZFC that $2^\omega = \kappa$ and for any subset $S \subset H$ of cardinality
$< \kappa$ the set of operators in $[B(H)]^+_1$ for which every nonzero
vector in the span of $S$ is cyclic is comeager in $[B(H)]^+_1$.}
\medskip

\noindent {\it Proof.} It is relatively consistent with ZFC + MA that
$2^\omega = \kappa$ ([6], Theorem 6.3). Let $S \subset H$ be a subset of
cardinality $< \kappa$. Then the span of $S$ is a union of fewer than
$\kappa$ finite dimensional subspaces of $H$, and the lemma implies that
for each such subspace $E$ the set of operators in $[B(H)]^+_1$ for which
every nonzero vector in $E$ is cyclic is comeager. The proof is completed
by observing that MA implies that the intersection of fewer than $2^\omega$
comeager sets in a Polish space is comeager ([3], Corollary 22C).\hfill\hal
\bigskip

Theorem 1 is related to the main theorem of [5]. That result implies, for
instance, that for any countable linearly independent subset $S$ of a
separable, infinite dimensional Hilbert space $H$ there exists $A \in B(H)$
for which every $v \in S$ is cyclic.
\medskip

Sophie Grivaux has pointed out to me that in Theorem 1 one can explicitly
construct an operator for which every nonzero vector in the span of $S$ is
cyclic. Namely, first find an orthonormal basis of $H$ whose span contains
$S$ (this can be accomplished by applying the Gram-Schmidt algorithm to $S$);
then one can show directly that the operator $V + V^*$, where $V$ is the
unilateral shift for the new basis, has the desired property. Moreover, with
a little work this idea can be used to show that the operators in $[B(H)]^+_1$
for which every nonzero vector in the span of $S$ is cyclic is dense for the
strong operator topology [4].
\medskip

I wish to thank Sophie Grivaux, Kenneth Kunen, and Peter Rosenthal for
helpful discussions.
\bigskip
\bigskip

[1] B.\ Beauzamy, {\it Introduction to Operator Theory and Invariant
Subspaces}, North-Holland (1988).
\medskip

[2] P.\ Enflo, On the invariant subspace problem for Banach spaces,
{\it Acta Math.\ \bf 158} (1987), 213-313. 
\medskip

[3] D.\ H.\ Fremlin, {\it Consequences of Martin's Axiom}, Cambridge
University Press (1984).
\medskip

[4] S.\ Grivaux, personal communication.
\medskip

[5] I.\ Halperin, C.\ Kitai, and P.\ Rosenthal, On orbits of linear
operators, {\it J.\ London Math.\ Soc.\ \bf 31} (1985), 561-565.
\medskip

[6] K.\ Kunen, {\it Set Theory: An Introduction to Independence Proofs},
Elsevier (1980).
\medskip

[7] C.\ J.\ Read, A solution to the invariant subspace problem on the space
$l_1$, {\it Bull.\ London Math.\ Soc.\ \bf 17} (1985), 305-317.
\bigskip
\bigskip

\noindent Math Dept.

\noindent Washington University

\noindent St.\ Louis, MO 63130 USA

\noindent nweaver@math.wustl.edu
\end